\documentclass[review]{elsarticle}

\usepackage{amsmath,amssymb}
\usepackage{hyperref}
\usepackage{mathrsfs}
\usepackage{amsmath,amscd}
\usepackage{fancyhdr}

\newtheorem{theorem}{Theorem}[section]
\newtheorem{proposition}[theorem]{Proposition}%
\newtheorem{lemma}[theorem]{Lemma}%

\begin{document}

\begin{frontmatter}

\title{Uniqueness and stability of normalized ground states for Hartree equation with a harmonic potential}
\author{Yi Jiang}
\address{School of Mathematical Sciences, Sichuan Normal University, Chengdu, 610066, China}

\author{Chenglin Wang}
\address{School of Science, Xihua University, Chengdu, 610039, China}

\author{Yibin Xiao}
\address{School of Mathematical Sciences, University of Electronic Science and Technology of China, Chengdu, 611731, China}


\author{Jian Zhang\corref{mycorrespondingauthor}}
\address{School of Mathematical Sciences, University of Electronic Science and Technology of China, Chengdu, 611731, China}
\cortext[mycorrespondingauthor]{Corresponding author}
\ead{zhangjian@uestc.edu.cn}

\author{Shihui Zhu\corref{}}
\address{School of Mathematical Sciences, Sichuan Normal University, Chengdu, 610066, China}

\begin{abstract}
The dynamic properties of normalized ground states for the Hartree equation with a harmonic potential are addressed. The existence of normalized ground state for any prescribed mass  is confirmed according to mass-energy  constrained variational approach. The uniqueness is shown by the  strictly convex properties of the energy functional. Moreover, the orbital stability of every normalized ground state is proven in terms  of the Cazenave and Lions' argument.
\end{abstract}


\end{frontmatter}


\section{Introduction}

In this paper, we study the following Hartree equation with a harmonic potential
\begin{equation}\label{1.1}
i\partial_{t}\psi+\Delta\psi-|x|^{2}\psi-(I_{\alpha}\ast|\psi|^{2})\psi=0,\;\;\; (t,x)\in \mathbb{R}\times\mathbb{R}^{N},
\end{equation}
where $N\geq1$ is the space dimension, and $I_{\alpha}:\mathbb{R}^{N}\backslash\{0\}\rightarrow\mathbb{R}$ is the Riesz potential defined by
\begin{equation}\label{1.2}
I_{\alpha}(x)=\frac{\Gamma(\frac{N-\alpha}{2})}{\Gamma(\frac{\alpha}{2})\pi^{\frac{N}{2}}2^{\alpha}|x|^{N-\alpha}}
\end{equation}
with $0<\alpha<N$ and $\Gamma$ is the Gamma function. Eq.(\ref{1.1}) arises typically if we consider the quantum mechanical time evolution of electrons in the mean field approximation of the many body effects, modeled by the Poisson equation, with a confinement modeled by the quadratic potential of the harmonic oscillator (see \cite{BBL1981,CMS2006,DK2003,GGKMP2003,K2004}).


\par From a mathematical point of view, Eq. (\ref{1.1}) without the term $|x|^{2}\psi$ is the standard defocusing Hartree equation, which does not possess any ground state,
but it  has a global solution with a scattering property in the natural energy space $H^{1}(\mathbb{R}^{N})$ (see \cite{C2003b,GV1980}). On the other hand, we recall the  focusing Hartree equation with a harmonic potential, which possesses ground states and blow-up solutions. Moreover, the  set stability of ground states and instability with blowing
-up of ground states are widely studied (see \cite{BDGS2014,CL1982,DS2014,GS1971,GV2012,MS2013,MS2017}). In fact, stability of every ground state is really concerned (see \cite{CL1982,F2001,FO2003b,F2004,JJTV2022,LR2020,O2018a,Z2005,Z2000a,Z2000b}). This  motivates us to further study the dynamics of  ground states for Eq.  (\ref{1.1}).

  The natural energy space to Eq. (\ref{1.1}) is defined by
\begin{equation}\label{1.3}
H(\mathbb{R}^{N}):=\Big\{u\in H^{1}(\mathbb{R}^{N}),\;\; \int_{\mathbb{R}^{N}}|x|^{2}|u|^{2}dx<\infty \Big\}.
\end{equation}
 For $u\in H(\mathbb{R}^{N})$, we define the mass functional
\begin{equation}\label{1.4}
	\mathcal{M}(u)=\int_{\mathbb{R}^{N}}|u|^{2}dx,
\end{equation}
and energy functional
\begin{equation}\label{1.5}
\mathcal{E}(u)=\frac{1}{2}\int_{\mathbb{R}^{N}}\big(|\nabla u|^{2}+|x|^{2}|u|^{2}+\frac{1}{2}(I_{\alpha}\ast|u|^{2})|u|^{2}\big)dx.
\end{equation}
Define the mass-energy constrained variational problem
\begin{equation}\label{1.6}
d(m):=\inf_{S}\mathcal{E}(u),
\end{equation}
where $S=\{u\in H(\mathbb{R}^{N}),\;\; \mathcal{M}(u)=m>0 \}$.

If the variational problem (\ref{1.6}) possesses a positive minimizer $Q$, then $Q$ must be a positive solution of the following Euler-Lagrange equation with the Lagrange multiplier $\omega$, which is a nonlinear elliptic equation
\begin{equation}\label{1.7}
-\Delta u+\omega u+|x|^{2}u+(I_{\alpha}\ast|u|^{2})u=0.
\end{equation}
We call $Q$ a ground state of Eq. (\ref{1.7}). Since
\begin{equation}\label{1.8}
\psi(t,x)=e^{i\omega t}Q(x)
\end{equation}
is a soliton solution of (\ref{1.1}), we also call $Q$ a ground state of (\ref{1.1}) and  $\omega$ the frequency of soliton.  Recently, the normalized solutions of nonlinear elliptic equations are widely studied  (see \cite{BMRV2021,BS2017,J1997,JL2022a,LZ2020,S2020a,S2020b,WW2022}).
A solution $u\in H(\mathbb{R}^{N})$ of Eq. (\ref{1.7}) is called the normalized solution  provided that $u$ satisfies the prescribed mass constraint $\mathcal{M}(u)=m$ for some $m>0$.  If the variational problem (\ref{1.6}) possesses a  minimizer $u$ for some $m>0$, then  $u$ must be a normalized solution of Eq. (\ref{1.7}) for some $m>0$.  A  positive minimizer   of  variational problem (\ref{1.6})   for some $m>0$ is called a normalized ground state   of  Eq. (\ref{1.7}) .  It is clear that a normalized ground state   of  Eq. (\ref{1.7})
is a normalized solution  of  Eq. (\ref{1.7}).

There are a plenty of known results for the existence of normalized solutions  and normalized ground states (see  \cite{BMRV2021,BS2017,GJ2018,J1997,JJTV2022,JL2021,JL2022a,LR2020,LZ2020,PR2025,S2020a,S2020b,WW2022}), but the uniqueness and dynamics
have  been less mentioned.  The main challenges come from the loss of scaling invariance of concerning nonlinear elliptic equations (see \cite{CL1982,JJTV2022,JL2022,LR2020,O2018a,Z2000b}).

In the present paper, we study  the convex properties of the energy functional $\mathcal{E}(u)$ and the set $S$ which is inspired by \cite{BC2013}. Then we imply  the  uniqueness of positive  minimizers of variational problem (\ref{1.6}). Thus
the uniqueness of normalized ground state $u$  for Eq. (\ref{1.7}) is got.  Moreover we further prove the orbital stability of every normalized ground state with respect to every frequency by applying Cazenave and Lions' argument (see \cite{CL1982}).

\par This paper is organized as follows. In Section 2, we state some preliminaries. In Section 3, we show the existence of normalized ground states. In Section 4, we  give the uniqueness of normalized ground states.  In Section 5, we   prove the orbital stability of every normalized ground state with respect to every frequency.

\section{Preliminaries}


We impose the initial data to Eq. (\ref{1.1}) as follows:
\begin{equation}\label{2.1}
	\psi(0,x)=\psi_{0}(x),\;\;\;x\in\mathbb{R}^{N}.
\end{equation}
According to Cazenave \cite{C2003b}, the  local well-posedness for the Cauchy problem (\ref{1.1}) and (\ref{2.1}) holds in $H(\mathbb{R}^{N})$.

\begin{proposition}\label{p2.1}
(\cite{C2003b}) Let $N\geq1$  and $0<\alpha<N$. For any $\psi_{0}(x)\in H(\mathbb{R}^{N})$, there exists $T>0$ such that the Cauchy problem (\ref{1.1}) and (\ref{2.1}) possesses a unique local solution
$\psi(t,x)\in \mathcal{C}([0,T),H(\mathbb{R}^{N}))$. In addition, the mass $\mathcal{M}(\psi)$
and the energy $\mathcal{E}(\psi)$ defined in (\ref{1.4}) and (\ref{1.5}), respectively, are conserved for all $t\in [0,T)$.
\end{proposition}

Next, we state the Hardy-Littlewood-Sobolev inequality (see \cite{LL2001}).

\begin{proposition}\label{p2.2}
(\cite{LL2001}) Let $\frac{1}{p}+\frac{1}{q}+\frac{\beta}{N}=2$ where $0<\beta<N$ and $1<p,q<\infty$. Then there exists a constant $C=C(N, p, \beta)$ such that for $f\in L^{p}(\mathbb{R}^{N})$ and $g\in L^{q}(\mathbb{R}^{N})$ one has $f(x)|x-y|^{-\beta}g(y)\in L^{1}(\mathbb{R}^{2N})$ and
\begin{equation*}
\Big|\int_{\mathbb{R}^{N}}\int_{\mathbb{R}^{N}}\frac{f(x)g(y)}{|x-y|^{\beta}}dxdy\Big|\leq C\|f\|_{L^{p}}\|g\|_{L^{q}}.
\end{equation*}
\end{proposition}

Then,  we state the following well-known convex properties (see \cite{Z}).

\begin{proposition}\label{p2.3}
(\cite{Z}) The functional $F:M\subseteq X\rightarrow \mathbb{R}$ has at most one minimum on $M$  in case the following hold:
\begin{itemize}
\item [(i)] $M$ is a convex subset of the linear space $X$.
\item [(ii)]  $F$ is strictly convex, i.e.
$$F((1-\lambda)u+\lambda v) <(1-\lambda) F(u)+\lambda F(v)$$
holds for all $u,v\in M$, $u\neq v$, and all $\lambda\in(0,1)$. \end{itemize}
\end{proposition}

Lastly, we give the following compactness result (see \cite{Z2000b}).
\begin{proposition}\label{p2.4}
(\cite{Z2000b}) Let $N\geq1$  and $2\leq q<\frac{2N}{(N-2)^+}$, where $\frac{2N}{(N-2)^{+}}=\infty$ for $N=1,2$, and $\frac{2N}{(N-2)^{+}}=\frac{2N}{N-2}$ for $N\geq 3$. Then, the embedding
\begin{equation}\label{2.2}
H(\mathbb{R}^{N})\hookrightarrow L^{q}(\mathbb{R}^{N})
\end{equation}
is compact.
\end{proposition}

\section{The existence  of the normalized ground state}

\begin{theorem}\label{L3.1}
Let $N\geq1$  and $0<\alpha<N$. Then, for any $m>0$, the
mass-energy variational problem (\ref{1.6})  is attained at  a non-negative  $Q\in S$.
\end{theorem}
Proof. Firstly the set $S$ is nonempty. Indeed, for any $u\in H(\mathbb{R}^{N})\backslash\{0\}$, we have $\mathcal{M}(u)=m_{0}$. Let
\begin{equation}\label{3.1}
v=\sqrt{\frac{m}{m_{0}}}u.
\end{equation}
It follows that $\mathcal{M}(v)=m$. Thus, $v\in S$.
\par  It is obvious that $d(m)\geq 0$ for any $u\in H(\mathbb{R}^{N})$. Let $\{u_{n}\}$ be a minimizing sequence of (\ref{1.6}); then, we have
\begin{equation}\label{3.2}
\mathcal{M}(u_{n})=m,\;\;\; \mathcal{E}(u_{n})\rightarrow d(m),\;\;\;\text{as}\;\;\;n\rightarrow\infty.
\end{equation}
By (\ref{1.5}), one has that $c>0$ such that
\begin{equation}\label{3.3}
\int_{\mathbb{R}^{N}}\big(|\nabla u_{n}|^{2}+|x|^{2}|u_{n}|^{2}+|u_{n}|^{2}\big)dx\leq d(m)+c,
\end{equation}
which shows that $\{u_{n}\}$ is bounded in $H(\mathbb{R}^{N})$. By the convexity estimate (see \cite{C2003b})
\begin{equation}\label{3.4}
\|\nabla |v|\|_{L^{2}(\mathbb{R}^{N})}\leq \|\nabla v\|_{L^{2}(\mathbb{R}^{N})}, \ \ v\in H^1(\mathbb{R}^{N}),
\end{equation}
 $\{|u_{n}|\}$ is also bounded in $H(\mathbb{R}^{N})$. Moreover, according to Proposition \ref{p2.4}, there exist  a subsequence $\{|u_{n}|\}$ (here, we use the same notation $\{|u_{n}|\}$) and a $H(\mathbb{R}^{N})$ function $u$ such that
\begin{equation}\label{3.5}
\begin{aligned}
&|u_{n}|\rightarrow Q \;\;\;\text{weakly}\;\;\; \text{in}\;\;\;H(\mathbb{R}^{N}),\\
&|u_{n}|\rightarrow Q \;\;\;\text{a.e.}\;\;\; \mathbb{R}^{N},\\
&|u_{n}|\rightarrow Q \;\;\;\text{strongly}\;\;\; \text{in}\;\;\;L^{p+1}(\mathbb{R}^{N}),
\end{aligned}
\end{equation}
where $1<p<\frac{N+2}{(N-2)^+}$. Since $0<\alpha<N$, it follows that
\begin{equation}\label{3.6}
|u_{n}|\rightarrow Q \;\;\;\text{strongly}\;\;\; \text{in}\;\;\;L^{\frac{4N}{N+\alpha}}(\mathbb{R}^{N}).
\end{equation}
From Proposition \ref{p2.2} and H\"{o}lder inequality, one has that
\begin{equation}\label{3.7}
\begin{aligned}
&\Big|\int_{\mathbb{R}^{N}}(I_{\alpha}\ast |u_{n}|^{2})|u_{n}|^{2}dx-\int_{\mathbb{R}^{N}}(I_{\alpha}\ast |Q|^{2})|Q|^{2}dx \Big|\\
=&\Big|\int_{\mathbb{R}^{N}}I_{\alpha}\ast (|u_{n}|^{2}+|Q|^{2})(|u_{n}|^{2}-|Q|^{2})dx \Big|\\
\leq& C\Big\||u_{n}|^{2}+|Q|^{2}\Big\|_{L^{\frac{2N}{N+\alpha}}}\Big\||u_{n}|^{2}-|Q|^{2}\Big\|_{L^{\frac{2N}{N+\alpha}}}\\
\leq &C\big(\|u_{n}\|_{L^{\frac{4N}{N+\alpha}}}^{3}+\|Q\|_{L^{\frac{4N}{N+\alpha}}}^{3}\big)\|u_{n}-Q\|_{L^{\frac{4N}{N+\alpha}}}.
\end{aligned}
\end{equation}
By (\ref{3.6}), it follows that
\begin{equation}\label{3.8}
\int_{\mathbb{R}^{N}}(I_{\alpha}\ast |u_{n}|^{2})|u_{n}|^{2}dx\rightarrow \int_{\mathbb{R}^{N}}(I_{\alpha}\ast |Q|^{2})|Q|^{2}dx\;\;\;\;\text{as}\;\;\;\;n\rightarrow\infty.
\end{equation}
We claim that $Q\neq 0$. In fact, if $Q=0$; then, by (\ref{3.5}),
\begin{equation}\label{3.9}
|u_{n}|\rightarrow 0, \;\;\;\text{in}\;\;\;L^{2}(\mathbb{R}^{N}).
\end{equation}
It follows that
\begin{equation}\label{3.10}
\lim_{n\rightarrow\infty}\mathcal{M}(|u_{n}|)=0,
\end{equation}
which contradicts $\mathcal{M}(|u_{n}|)=m>0$. Thus, $Q\neq 0$.
\par Using the weakly lower semicontinuous property of norms, it follows that
\begin{equation}\label{3.11}
\liminf_{n\rightarrow\infty}\|\nabla |u_{n}|\|_{L^{2}}^{2}\geq \|\nabla Q\|_{L^{2}}^{2},
\end{equation}
\begin{equation}\label{3.12}
\liminf_{n\rightarrow\infty}\int_{\mathbb{R}^{N}}|x|^{2}|u_{n}|^{2}dx\geq \int_{\mathbb{R}^{N}}|x|^{2}|Q|^{2}dx.
\end{equation}
Thus, one has
\begin{equation}\label{3.14}
d(m)=\liminf_{n\rightarrow\infty}\mathcal{E}(|u_{n}|)\geq \mathcal{E}(Q).
\end{equation}
Since
\begin{equation}\label{3.15}
|u_{n}|\rightarrow Q\;\;\;\text{in}\;\;\; L^{2}(\mathbb{R}^{N}),
\end{equation}
then it yields that $\mathcal{M}(Q)=m$. By the definition of $d(m)$, it is verified that
\begin{equation}\label{3.16}
d(m)\leq \mathcal{E}(Q)\leq \liminf_{n\rightarrow\infty}\mathcal{E}(|u_{n}|)=d(m).
\end{equation}
This implies that $Q$ is a non-negative minimizer of the variational problem (\ref{1.6}).

\begin{theorem}\label{3}
Let $N\geq1$  and $0<\alpha<N$. Then, for any $m>0$,  Eq. (\ref{1.7}) possesses a normalized ground state  for the prescribed mass constraint $\mathcal{M}(u)=m$. Moreover, one has that the Lagrange multiplier $\omega<0$.
\end{theorem}
Proof. In terms of  Lemma 3.1, the variational problem (\ref{1.6}) possesses a non-negative minimizer $Q$ such that  $\mathcal{M}(Q)=m$. Then, $Q$ satisfies the following Euler-Lagrange equation
\begin{equation}\label{3.2a}
-\Delta u+\omega u+|x|^{2}u+(I_{\alpha}\ast|u|^{2})u=0
\end{equation}
for some $\omega\in \mathbb{R}$. Applying the strong  maximum principle, $Q\geq0 $ and $\mathcal{M}(Q)=m>0$ imply that $Q>0$. Therefore,  for any $m>0$,  Eq. (\ref{1.7}) possesses a normalized ground state  for the prescribed mass constraint $\mathcal{M}(u)=m$. By (\ref{3.2a}), it follows that
\begin{equation}\label{3.2b}
\int | \nabla u|^2dx+\omega \int | u|^2dx+\int |x|^{2}|u|^2dx+\int (I_{\alpha}\ast|u|^{2})|u|^{2}dx=0,
\end{equation}
which implies that $\omega<0$.

\section{Uniqueness of the normalized ground state}

\begin{lemma}\label{L41}
Let $\widetilde{S}$ be defined as
\begin{equation*}
\widetilde{S}=\big\{\rho\in L^{1}(\mathbb{R}^{N}):  \; \rho\geq 0,\ \  \int_{\mathbb{R}^{N}}\rho dx=m>0\big\}.
\end{equation*}
Then, the set $\widetilde{S}$ is convex in  $L^{1}(\mathbb{R}^{N})$.
\end{lemma}
Proof.  For any
\begin{equation}\label{4.1a}
\rho_{1}\in \widetilde{S}, \;\;\; \rho_{2}\in\widetilde{S},
\end{equation}
and for all $0<\lambda<1$, one can check that $\lambda \rho_{1}+(1-\lambda)\rho_{2}\in L^{1}(\mathbb{R}^{N})$, and
 it follows that
\begin{equation}\label{4.2b}
\int_{\mathbb{R}^{N}}\lambda \rho_{1}+(1-\lambda)\rho_{2}dx=m.
\end{equation}
By (\ref{4.2b}) and the definition of $\widetilde{S}$, one deduces that $\lambda\rho_{1}+(1-\lambda)\rho_{2}\in \widetilde{S}$. Therefore, $\widetilde{S}$ is a convex set.

\begin{lemma}\label{L42}
Let  $\widetilde{\mathcal{E}}(\sqrt{\rho})$ be defined as
\begin{equation*}
\widetilde{\mathcal{E}}(\sqrt{\rho})=\frac{1}{2}\int_{\mathbb{R}^{N}}\big(|\nabla \sqrt{\rho}|^{2}+|x|^{2}\rho+\frac{1}{2}(I_{\alpha}\ast\rho)\rho\big)dx.
\end{equation*}
Then, the functional  $\widetilde{\mathcal{E}}(\rho)$ is strictly convex with respect to $\rho$.
\end{lemma}
Proof.   From (\ref{4.1a}), without loss of generality, one can assume that
\begin{equation}\label{4.2a}
u_{1}=\sqrt{\rho_{1}}\;\;\;\text{and}\;\;\; u_{2}=\sqrt{\rho_{2}},
\end{equation}
and for $0<\lambda<1$, denote
\begin{equation}\label{4.2b}
\lambda \rho_{1}+(1-\lambda)\rho_{2}=|u_{\lambda}|^{2}.
\end{equation}
It follows that
\begin{equation}\label{4.2c}
u_{\lambda}=\sqrt{\lambda \rho_{1}+(1-\lambda)\rho_{2}}=\sqrt{\lambda u_{1}^{2}+(1-\lambda)u_{2}^{2}}\in S.
\end{equation}
Applying Cauchy inequality, one deduces that
\begin{equation}\label{4.2d}
\begin{aligned}
u_{\lambda}\nabla u_{\lambda}=&\lambda u_{1}\nabla u_{1}+(1-\lambda)u_{2}\nabla u_{2}\\
=&(\sqrt{\lambda}u_{1})(\sqrt{\lambda}\nabla u_{1})+(\sqrt{1-\lambda}u_{2})(\sqrt{1-\lambda}\nabla u_{2})\\
\leq &\sqrt{\lambda u_{1}^{2}+(1-\lambda)u_{2}^{2}} \sqrt{\lambda |\nabla u_{1}|^{2}+(1-\lambda)|\nabla u_{2}|^{2}}\\
=&u_{\lambda}\sqrt{\lambda |\nabla u_{1}|^{2}+(1-\lambda)|\nabla u_{2}|^{2}}.
\end{aligned}
\end{equation}
Since $u_{\lambda}>0$, it follows that
\begin{equation}\label{4.2e}
\nabla u_{\lambda}\leq \sqrt{\lambda |\nabla u_{1}|^{2}+(1-\lambda)|\nabla u_{2}|^{2}}.
\end{equation}
From (\ref{4.2a}) and (\ref{4.2c}), one has that
\begin{equation}\label{4.2f}
\big|\nabla \sqrt{\lambda \rho_{1}+(1-\lambda)\rho_{2}}\big|^{2} \leq \lambda |\nabla \sqrt{\rho_{1}}|^{2}+(1-\lambda)|\nabla \sqrt{\rho_{2}}|^{2},
\end{equation}
which shows that $\int_{\mathbb{R}^{N}}|\nabla\sqrt{\rho}|^{2}dx$ is strictly convex in $\rho$.
\par Next we claim that $\int_{\mathbb{R}^{N}}(I_{\alpha}\ast\rho)\rho dx$ is a convex in $\rho$. In fact, for $0<\lambda<1$ one has that
\begin{equation}\label{4.2g}
\begin{aligned}
&\int_{\mathbb{R}^{N}}I_{\alpha}\ast\big(\lambda\rho_{1}+(1-\lambda)\rho_{2}\big))[\lambda\rho_{1}+(1-\lambda)\rho_{2}]dx\\
=&\int_{\mathbb{R}^{N}}\int_{\mathbb{R}^{N}}\frac{\lambda\rho_{1}(x)+(1-\lambda)\rho_{2}(x)}{|x-y|^{N-\alpha}}[\lambda\rho_{1}(y)+(1-\lambda)\rho_{2}(y)]dxdy\\
=&\int_{\mathbb{R}^{N}}\int_{\mathbb{R}^{N}}\frac{\lambda^{2}\rho_{1}(x)\rho_{1}(y)+2\lambda(1-\lambda)\rho_{1}(x)\rho_{2}(y)+(1-\lambda)^{2}\rho_{2}(x)\rho_{2}(y)}{|x-y|^{N-\alpha}}dxdy\\
\leq & \int_{\mathbb{R}^{N}}\int_{\mathbb{R}^{N}}\frac{\lambda\rho_{1}(x)\rho_{1}(y)+(1-\lambda)\rho_{2}(x)\rho_{2}(y)}{|x-y|^{N-\alpha}}dxdy\\
= & \lambda \int_{\mathbb{R}^{N}}\int_{\mathbb{R}^{N}}\frac{\rho_{1}(x)\rho_{1}(y)}{|x-y|^{N-\alpha}}dxdy+(1-\lambda)\int_{\mathbb{R}^{N}}\int_{\mathbb{R}^{N}}\frac{\rho_{2}(x)\rho_{2}(y)}{|x-y|^{N-\alpha}}dxdy\\
=&\lambda\int_{\mathbb{R}^{N}}(I_{\alpha}\ast \rho_{1})\rho_{1}dx+(1-\lambda)\int_{\mathbb{R}^{N}}(I_{\alpha}\ast \rho_{2})\rho_{2}dx.
\end{aligned}
\end{equation}
\par Since $\int_{\mathbb{R}^{N}}|x|^{2}\rho dx$ is strictly convex in $\rho$, $\int_{\mathbb{R}^{N}}|\nabla\sqrt{\rho}|^{2}dx$ is strictly convex in $\rho$ and $\int_{\mathbb{R}^{N}}(I_{\alpha}\ast\rho)\rho dx$ is a strictly convex in $\rho$, then one gets that $\widetilde{\mathcal{E}}(\sqrt{\rho})$  is strictly convex in $\rho$.

\begin{theorem}\label{t3.1}
Let $N\geq1$  and $0<\alpha<N$. Then, for any $m>0$,  Eq. (\ref{1.7})  possesses a unique  normalized ground state  with the prescribed mass $\mathcal{M}(u)=m$. Moreover, all minimizers of the variational problem (\ref{1.6}) are in the set $\{Qe^{i\theta}:\ \ \theta\in \mathbb{R}\}$.
\end{theorem}
Proof.  According to Theorem  \ref{3}, the mass-energy constrained variational problem (\ref{1.6}):
\begin{equation*}
d(m):=\inf_{u\in S}\mathcal{E}(u),
\end{equation*}
where $S=\{u\in H(\mathbb{R}^{N}),\;\; \int |u|^2dx=m>0 \}$, possesses a normalized ground state  for the prescribed mass constraint $\mathcal{M}(u)=m>0$.

We take the following transformation:
\begin{equation}\label{4.3a}
\rho(x)=|u(x)|^{2}.
\end{equation}
Then, we can see that $\rho\in \widetilde{S}$ implies that $ u\in S$ and $\widetilde{\mathcal{E}}(\sqrt{\rho})=\mathcal{E}(u)$. The mass-energy constrained variational problem (\ref{1.6}) can be rewritten as following form
\begin{equation}\label{4.3b}
\widetilde{d}(\rho)=\inf_{\rho\in \widetilde{S}}\widetilde{\mathcal{E}}(\sqrt{\rho}).
\end{equation}
It follows from  Lemma \ref{L41} and Lemma \ref{L42} that the variational problem (\ref{4.3b}) is the strictly convex functional  $\widetilde{\mathcal{E}}(\sqrt{\rho})$ defined on the convex set $\widetilde{S}$. Applying Proposition \ref{p2.3}, one deduces that the  variational problem (\ref{4.3b}) possesses at most one minimum point in $\widetilde{S}$. Then, from (\ref{4.3a})  the variational problem (\ref{1.6}) possesses at most one positive minimizer. Therefore, Eq. (\ref{1.7})  possesses a unique  normalized ground state  with the prescribed mass $\mathcal{M}(u)=m$.

Suppose that $v$ is  a minimizer of the variational problem (\ref{1.6}). One has that
\begin{equation}\label{4.3e}
v=|v|e^{i\theta}\;\;\;\text{for some}\;\;\; \theta\in \mathbb{R}.
\end{equation}
Since for $v\in H(\mathbb{R}^{N})$,
\begin{equation}\label{4.3f}
\int\rvert\nabla v\rvert^{2}dx\geq \int\rvert\nabla \rvert v\rvert\rvert^{2}dx,
\end{equation}
it follows that
\begin{equation}\label{4.3g}
E(v)\geq E(\rvert v\rvert).
\end{equation}
It yields that $|v|$ is also a minimizer of the variational problem (\ref{1.6}). Then one implies that
\begin{equation}\label{4.3h}
|v|=Q.
\end{equation}
It follows that
\begin{equation}\label{4.3i}
v\in \{Qe^{i\theta}:\ \ \theta\in \mathbb{R}\}.
\end{equation}

\section{Orbital stability of the normalized ground state}

\begin{lemma}\label{L5.1}
Let $N\geq1$  and $0<\alpha<N$. Suppose that the initial data $\psi_{0}(x)\in H(\mathbb{R}^{N})$. Then, the corresponding solution $\psi(t,x)$ to   Cauchy problem (\ref{1.1}) and (\ref{2.1})  exists globally for all time $t$.
\end{lemma}
Proof. Since  the initial $\psi_{0}(x)\in H(\mathbb{R}^{N})$, from Proposition \ref{p2.1},
 there exists $T>0$ such that the Cauchy problem (\ref{1.1}) and (\ref{2.1}) possesses a unique local solution
$
\psi(t,x)\in \mathcal{C}([0,T),H(\mathbb{R}^{N}))
$, which satisfies two conservation laws. Then, for all time $t$,
 \begin{equation}\label{5.1a}
 \begin{array}{lll}\vspace{0.3cm}
 \mathcal{E}(\psi_0)&=\mathcal{E}(\psi(t))\\
 \vspace{0.3cm}
 &=\frac{1}{2}\int_{\mathbb{R}^{N}}\big(|\nabla \psi(t)|^{2}+|x|^{2}|\psi(t)|^{2}+\frac{1}{2}(I_{\alpha}\ast|\psi(t)|^{2})|\psi(t)|^{2}\big)dx\\

 &\geq \frac 12 \|\nabla \psi(t)\|_{L^2(\mathbb{R}^{N})}^2+\frac{1}{2}\||x|\psi(t)\|_{L^2(\mathbb{R}^{N})}^2.
\end{array}
\end{equation}
By combining above estimates with the conservation of mass, there exists a positive constant $C>0$ such that
\begin{equation}\label{5.1b}
\|\nabla \psi(t)\|_{L^2(\mathbb{R}^{N})}^2+\||x|\psi(t)\|_{L^2(\mathbb{R}^{N})}^2+\| \psi(t)\|_{L^2(\mathbb{R^{N})}}^2\leq C
\end{equation}
holds for all time, which implies that the solution  $\psi(t)$ to   Cauchy problem (\ref{1.1}) and (\ref{2.1})  exists globally for all time $t$.

\begin{theorem}\label{t5.2}
Let $N\geq1$  and $0<\alpha<N$. Then, for all  $m>0$,  Eq. (\ref{1.7})  possesses a unique  normalized ground state $Q$ with the prescribed mass $\mathcal{M}(Q)=m$. Moreover, $Q$ is  orbitally stable under the evolution flow of Eq. (\ref{1.1}). That is,   for arbitrary
$\varepsilon>0$, there exists $\delta>0$ such that to the initial data $\psi_0\in H(\mathbb{R}^{N})$, if
\begin{equation}\label{5.3a}
\inf\limits_{\theta\in\mathbb{R}}\|\psi_0(\cdot)-Q(\cdot)e^{i\theta}\|_{H(\mathbb{R}^{N})}<\delta,
\end{equation}
then the corresponding solution $\psi(t,x)$ of Cauchy
problem (\ref{1.1}) satisfies
\begin{equation}\label{5.3b}
\inf\limits_{\theta\in\mathbb{R}}\|\psi(t,\cdot)-Q(\cdot)e^{i\theta}\|_{H(\mathbb{R}^{N})}<\varepsilon  \end{equation}
for all $t>0$.
\end{theorem}
Proof.  Since $\psi_0\in H(\mathbb{R}^{N})$, from Lemma\  \ref{L5.1},
the corresponding  solution $\psi(t,x)$ of  Cauchy problem
(\ref{1.1}) exists globally  in $H(\mathbb{R}^{N})$.

Now, we prove Theorem\  \ref{t5.2} \ by contradiction in terms of
Cazenave and Lions' arguments (see also
\cite{C2003b}). Assume the conclusion in Theorem  \ref{t5.2}  does not hold. That is,
  there
exist  $\varepsilon_0>0$ and a series of initial sequence
$\{\psi_{0,n}\}_{n=1}^{+\infty}$ such that
\begin{equation}\label{5.2c}
\inf\limits_{\theta\in
\mathbb{R}} \|\psi_{0,n}(x)-Q(x)e^{i\theta}\|_{H(\mathbb{R}^{N})}<\frac1n,
\end{equation}
and there exists $\{t_n\}_{n=1}^{+\infty}$  such that the
corresponding solution sequence  $\{\psi_n(t_n,x)\}_{n=1}^{+\infty}$ of the Cauchy problem (\ref{1.1})
satisfies
\begin{equation}\label{5.2d}
\inf\limits_{\theta\in\mathbb{R}} \|\psi_n(t_n,x)-Q(x)e^{i\theta}\|_{H(\mathbb{R}^{N})}\geq\varepsilon_0
\end{equation}
for any $\theta\in \mathbb{R}$. From    (\ref{5.2c}) and the conservation
laws in Proposition \ref{p2.1}, we see that  as $n\rightarrow+\infty$,
\[\int |\psi_n(t_n,x)|^2dx=\int |\psi_{0,n}|^2dx\rightarrow \int |v|^2dx=m,\]
\[\mathcal{E}(\psi_n(t_n,x))=\mathcal{E}(\psi_{0,n})\rightarrow \mathcal{E}(v)=d(m).\]
Hence,  $\{\psi_n(t_n,x)\}_{n=1}^{+\infty}$ is a minimizing
sequence of the   variational problem (\ref{1.6}).   Therefore, by Theorem \ref{t3.1} there exists a $v\in \{Q(x)e^{i\theta}:\ \ \theta\in \mathbb{R}\}$  such that
\begin{equation}\label{5.3e}
\|\psi_n(t_n,x)-v(x)\|_{H(\mathbb{R}^{N})}\rightarrow 0\  \ {\rm as}\ \ n\rightarrow+\infty. \end{equation}
We see that   (\ref{5.3e})  contradicts with (\ref{5.3b}). Then, the conclusion in Theorem\  \ref{t5.2} \  is true.  This completes the proof.\\


\section*{Acknowledgment}

This research is supported by the National Natural Science Foundation of China 12271080, 12571318 and Sichuan Technology Program 25LHJJ0156.

\section*{Declarations}

All of authors state no conflict of interest and no data was used for the research described in the article.



\begin{thebibliography}{10}

	
	
	
\bibitem{BC2013}W. Bao and Y. Cai, Mathematical theory and numerical methods for Bose-Eisntein condensation, \emph{Kinet. Relat. Mod.} \textbf{61}(2013), 1-135. 
\bibitem{BDGS2014}D. Bonanno, P. d'Avenia, M. Ghimenti and M. Squassina, Soliton dynamics for the generalized Choquard equation, \emph{J. Math. Anal. Appl.} \textbf{417}(2014), 180-199.
\bibitem{BMRV2021}T. Bartsch, R. Molle, M. Rizzi and G. Verzini, Normalized solutions of mass supercritical Schr\"odinger equations with potential, \emph{Commun. Partial Differ. Equ.} \textbf{46}(2021), 1729-1756.
\bibitem{BS2017}T. Bartsch and N. Soave, A natural constraint approach to normalized solutions of nonlinear Schr\"odinger equations and systems, \emph{J. Funct. Anal.} \textbf{272}(2017), 4998-5037.
\bibitem{BBL1981}R. Benguria, H. Brezis and E. H. Lieb, The Thomas-Fermi-von Weizs\"{a}cker theory of atoms and molecules, \emph{Commun. Math. Phys.} \textbf{79}(1981), 167-180.
\bibitem{CMS2006} R. Carles, N. J. Mauser and H. P. Stimming, (Semi) Classical limit of the Hartree equation with harmonic potential, \emph{SIAM J. Appl. Math.} \textbf{3}(2006), 2112758. 
\bibitem{C2003b}T. Cazenave, Semilinear Schr\"{o}dinger equations, Courant Institute of Mathematical Sciences, American Mathematical Society, New York, (2003).

\bibitem{CL1982} T. Cazenave,  P. L. Lions, Orbital stability of standing waves for some nonlinear Schr\"odinger equations. \emph{Commun. Math. Phys.} \textbf{85}(1982), 549-561.

\bibitem{DK2003} B. Deconinck, J. N. Kutz, Singular instability of exact stationary solutions of the non-local Gross-Pitaevskii equations. \emph{Phys. Lett. A} \textbf{319}(2003), 97-103.
\bibitem{DS2014} P. d'Avenia and M. Squassina, Soliton dynamics for the Schr\"{o}dinger-Newton system, \emph{Math. Models Methods Appl. Sci.} \textbf{24}(2014), 553-572.
\bibitem{F2001}R. Fukuizumi, Stability and instability of standing waves for the nonlinear Schr\"{o}dinger equation with harmonic potential, \emph{Discrete Contin. Dyn. Syst.} \textbf{7}(2001), 525-544. 

\bibitem{FO2003b}R. Fukuizumi and M. Ohta, Stability of standing waves for nonlinear Schr\"{o}dinger equations with potentials, \emph{Differ. Integral. Equ.} \textbf{16}(2003), 111-128.
\bibitem{F2004}R. Fukuizumi, Stability of standing waves for nonlinear Schr\"{o}dinger equations with critical power nonlinearity and potentials, \emph{Adv. Differ. Equ.} \textbf{10}(2004), 259-276.
\bibitem{GGKMP2003}J. J. Garcia-Ripoll, V. V. Konotop, B. Malomed and V. M. Perez-Garcia, A quasi-local Gross-Pitaevskii for Bose-Einstein condensates, \emph{Math. Compt. Simulation}, \textbf{62}(2003), 21-30. 
\bibitem{GS1971}K. Gustafson and D. Sather, A branching analysis of the Hartree equation, \emph{Rend. Mat.} \textbf{4}(1971), 21-30.
\bibitem{GV1980}J. Ginibre and G. Velo, On a class of nonlinear Schr\"{o}dinger equations, \emph{Math. Z.} \textbf{170}(1980) 109-136.
\bibitem{GV2012}H. Genev and G. Venkov, Soliton and blow-up solutions to the time-dependent Schr\"{o}dinger-Hartree equation, \emph{Disc. Contin. Dyn. Syst. Ser. S} \textbf{5}(2012) 903-923.
\bibitem{GJ2018}T. Gou and L. Jeanjean, Multiple positive normalized solutions for nonlinear Schr\"odinger systems, \emph{Nonlinearity}, \textbf{31}(2018), 2319-2346. 

\bibitem{J1997}L. Jeanjean, Existence of solutions with prescribed norm for semilinear elliptic equations, \emph{Nonlinear Anal., Theory Methods Appl.} \textbf{28}(1997), 1633-1659. 
\bibitem{JJTV2022}L. Jeanjean, J. Jendrej, T. Le and N. Visciglia, Orbital stability of ground states for a Sobolev critical Schr\"{o}dinger equation,\emph{ J. Math. Pures Appl.} \textbf{164}(2022), 158-179.
\bibitem{JL2021}L. Jeanjean and T. Le, Multiple normalized solutions for a Sobolev critical Schr\"odinger equation, \emph{Math. Ann. }\textbf{384}(2022), 101-134.
\bibitem{JL2022a}L. Jeanjean and S. S. Lu, On global minimizers for a mass constrained problem, \emph{Calc. Var. Partial Differ. Equ.} \textbf{61}(2022), 214.
\bibitem{JL2022}L. Jeanjean, S. S. Lu. Normalized solutions with positive energies for a coercive problem and application to the cubic-quintic nonlinear Schr\"{o}dinger equation. \emph{Math. Mod. Meth. Appl. Sci.} \textbf{32}(2022) 1557-1588.

\bibitem{K2004}M. Kurth, On the existence of infinitely many modes of a nonlocal nonlinear Schr\"{o}dinger equation related to Dispersion-Managed solitons, \emph{SIAM J. Math. Anal.} \textbf{36}(2004), 967-985.
\bibitem{LR2020}M. Lewin and S. Rota Nodari, The double-power nonlinear Schr\"{o}dinger equation and its generalizations: uniqueness, non-degeneracy and applications, \emph{Calc. Var. Partial Differ. Equ.} \textbf{59}(2020), 197.
\bibitem{LL2001}E. H. Lieb and M. Loss, Analysis, Second edition, Graduate Studies in Mathematics, 14, American Mathematical Society, Providence, RI, 2001.

\bibitem{LZ2020} H. Luo and Z. Zhang, Normalized solutions to the fractional Schr\"{o}dinger equations with combined nonlinearities, \emph{Calc. Var. Partial Differ. Equ.} \textbf{59}(2020), Paper no.143, 35pp.
\bibitem{MS2013} V. Moroz and J. Van Schaftingen, Ground states of nonlinear Choquard equations: existence, qualitative properties and decay asymptotics, \emph{J. Funct. Anal.} \textbf{265}(2013), 153-184.
\bibitem{MS2017} V. Moroz and J. Van Schaftingen, A guide to the Choquard equation, \emph{J. Fixed Point Theory. Appl.} \textbf{19}(2017), 773-813.

\bibitem{O2018a}M. Ohta, Strong instability of standing waves for nonlinear Schr\"{o}dinger equations with harmonic potential, \emph{Funkc. Ekvacioj, } \textbf{61}(2018), 135-143. 
\bibitem{PR2025}X. Peng and M. Rizzi, Normalized solutions of mass supercritical Schr\"{o}dinger-Posson equation with potential, \emph{Calc. Vari. Partial Differ. Equ.} \textbf{64}(2025), no. 152.

\bibitem{S2020a}N. Soave, Normalized ground states for the NLS equation with combined nonlinearities, \emph{J. Differ. Equ.} \textbf{269}(2020), 6941-6987.
\bibitem{S2020b}N. Soave, Normalized ground states for the NLS equation with combined nonlinearities: The Sobolev critical case, \emph{J. Funct. Anal.} \textbf{279}(2020), 108610.
\bibitem{WW2022}J. C. Wei and Y. Z. Wu, Normalized solutions for Schr\"{o}dinger equations with critical Sobolev exponent and mixed nonlinearities, \emph{J. Funct. Anal.} \textbf{283}(2022), 109574.
\bibitem{Z2005} J. Zhang, Sharp threshold for blowup and global existence in nonlinear Schr\"{o}dinger equations under a harmonic potential, \emph{Commun Partial Differ. Equ.} \textbf{30}(2005), 1429-1443. 
\bibitem{Z2000a}J. Zhang, Stability of attractive Bose-Einstein condenstates, \emph{J. Stat. Phys.} \textbf{101}(2000), 731-746. 
\bibitem{Z2000b}J. Zhang, Stability of standing waves for nonlinear Schr\"{o}dinger equations with unbounded potentials, \emph{Z. Angew. Math. Phys.} \textbf{51}(2000), 498-503. 
\bibitem{Z} E. I. N.  Zeidler, Nonlinear Functional Analysis and Its Applications,  Springer New York (1988).


\end{thebibliography}
\end{document}